\documentclass[11pt]{article}
\usepackage{latexsym}
\usepackage{amsmath}
\usepackage{amssymb}
\usepackage{amsfonts}

\numberwithin{equation}{section}

\newcommand{\beq}{\begin{equation}}
\newcommand{\eeq}{\end{equation}}
\def\nm{\noalign{\medskip}}

\newcommand{\Cbb}{{\mathbb C}}

\newcommand{\Rbb}{{\mathbb R}}

\newcommand{\p}{\partial}

\newcommand{\eps}{\varepsilon}
\newcommand{\R}{{\mathbb R}}

\renewcommand{\div}{{\rm div}}
\newcommand{\la}{\langle}
\newcommand{\ra}{\rangle}

\newcommand{\Bn}{{\bf n}}
\newcommand{\Bt}{{\bf t}}
\newcommand{\Bv}{{\bf v}}
\newcommand{\Bx}{{\bf x}}
\newcommand{\BA}{{\bf A}}

\newcommand{\BI}{{\bf I}}

\newcommand{\BR}{{\bf R}}

\newcommand{\Gs}{\sigma}

\newcommand{\GO}{\Omega}

\title{Equivalence of inverse problems for 2D elasticity and for the thin plate with finite measurements and its applications}
\author{Hyeonbae Kang\thanks{Department of Mathematics, Inha University, Incheon 402-751, Korea (hbkang@inha.ac.kr).}\qquad Graeme Milton\thanks{Department of Mathematics, University of Utah, 155 S 1400 E RM 233, Salt Lake City, Utah 84112, USA (milton@math.utah.edu).}\qquad Jenn-Nan Wang\thanks{Department of Mathematics, NCTS (Taipei), National Taiwan University, Taipei 106, Taiwan (jnwang@math.ntu.edu.tw).}}

\date{}

\begin{document}
\maketitle

\begin{abstract}
In this paper, we prove that the inverse problems for 2D elasticity and for the thin plate with boundary data (finite or full measurements) are equivalent. Having proved this equivalence, we can solve inverse problems for the plate equation with boundary data by solving the corresponding inverse problems for 2D elasticity, and vice versa. For example, we can derive bounds on the volume fraction of the
two-phase thin plate from the knowledge of one pair of boundary measurements using the known result for 2D elasticity \cite{ml}. Similarly, we give another approach to the size estimate problem for the thin plate studied by  Morassi, Rosset, and Vessella \cite{mrv07, mrv09}.
\end{abstract}

\section{Introduction}\label{sec1}
\setcounter{equation}{0}

In this note we would like to connect the inverse boundary value problem for thin plate to that for 2D elasticity. To begin, we first discuss the inverse boundary value problem for the thin plate.  Let $\Omega\subset\R^2$ be a simply-connected bounded domain with smooth boundary
$\partial\Omega$. Supposed that the middle surface of the thin plate
with uniform thickness $h$ occupies $\Omega$. In the Kirchhoff-Love
theory of thin elastic plates, the
transversal displacement $u$ satisfies
\begin{equation}\label{plate}
\text{div}\ \text{div}({\mathbb
C}\nabla^2u)=0\quad\text{in }\Omega,
\end{equation}
where $\nabla^2u$ is the Hessian matrix of $u$, i.e.,
$$\nabla^2u=\begin{bmatrix}u_{,11}&u_{,12}\\u_{,12}&u_{,22}\end{bmatrix},$$
and ${\mathbb C}=(C_{ijkl}(\Bx)), i,j,k,l=1,2,$ is a 4th order tensor satisfying that ${\mathbb C}\in L^{\infty}(\Omega)$,
 \beq\label{sym}
\text{(symmetry property)}\quad C_{ijkl}(\Bx)=C_{jikl}(\Bx)=C_{klij}(\Bx)\quad\forall\ \Bx\in\Omega\ a.e.,
 \eeq
 and there exists $\gamma>0$ such that
 \beq\label{con}
\text{(strong convexity)}\quad {\mathbb C}\BA\cdot \BA\ge\gamma |\BA|^2\quad\forall\ \Bx\in\Omega\ a.e.,
 \eeq
for every $2\times 2$ symmetric matrix $\BA$. Hereafter, for any function $u$, $u_{,j}$ denotes the derivative of $u$ with respect to $x_j$.  The Dirichlet data associated with \eqref{plate} is described by the pair $\{u,u_{,n}\}$ and the
Neumann data by the pair
\[
 ({\mathbb C}\nabla^2u){\bf n}\cdot{\Bn}=-M_n, \quad \text{div}({\mathbb C}\nabla^2u)\cdot{\Bn} + (({\mathbb C}\nabla^2u){\bf n}\cdot{\Bt})_{,t}=(M_{t})_{,t},
 \]
where ${\bf
n}$ is the boundary normal, ${\bf t}$ is the unit tangent
vector field along $\partial\Omega$ in the positive orientation, $u_{,n}=\nabla u\cdot{\bf n}$, and
$u_{,t}=\nabla u\cdot{\bf t}$. The quantities $M_n$ and
$M_t$ are known as the twisting moment and the bending moment
applied on the boundary $\partial\Omega$.

The 2D elasticity equation is given by
\begin{equation}\label{elasticity}
\text{div}\ \sigma=0\quad\text{in}\quad\Omega,
\end{equation}
where $\sigma$ is the stress, which is related to the strain $\eps=(\nabla \Bv+(\nabla \Bv)^T)/2$ ($T$ for transpose) by Hooke's law
\[
\eps={\mathbb S}\sigma.
\]
Here $\Bv$ is the displacement field and ${\mathbb S}$ is known as the compliance tensor. We assume that ${\mathbb S}\in L^{\infty}(\Omega)$, and that ${\mathbb S}$ satisfies the symmetry condition \eqref{sym} and the strong convexity condition \eqref{con} (with a possibly different constant). The Dirichlet and Neumann conditions for \eqref{elasticity} are described by
$\Bv$ and $\sigma{\bf n}$, respectively.

It is widely known that \eqref{plate} and \eqref{elasticity} are
equivalent (see for example \cite{mbook}) as long as $\Omega$ is
simply connected. As for the equivalence of boundary data, Ikehata
\cite{ikehata} showed that when the full measurements (the
Dirichlet-to-Neumann map or the Neumann-to-Dirichlet map) are
allowed, knowing one set of boundary data uniquely determines the
other one.  Consequently, the two inverse boundary value problems
with the full measurements are equivalent. Note that Ikehata's
result was a uniqueness proof. We would also like to mention that
explicit formulas for constructing the Dirichlet data $(u, u_{,n})$ of the thin
plate from the traction data $\Gs\Bn$ of 2D elasticity were given in
\cite[p.158]{gu}. However, to the best of our knowledge, it does not
seem that the relation between the individual boundary data $\{\Bv\}$ and $\{M_n, M_t \}$ is known.

It is the purpose of this paper to clarify the relation between the boundary data $\{\Bv, \Gs\Bn \}$ and $\{ u, u_{,n}, M_n, M_t \}$. In fact, we show that $\Bv$ on $\partial\GO$ explicitly determines $(M_n, M_t)$, and, following \cite{gu}, show that
$\Gs\Bn$ on $\partial\Omega$ explicitly determines $(u, u_{,n})$, and vice versa. As a consequence, the Dirichlet (resp. Neumann) boundary value problem for 2D elasticity equation is equivalent to the Neumann (resp. Dirichlet) boundary value problem for the thin plate equation, and thus the two inverse boundary value problems with finite measurements are equivalent. We emphasize that the results of this paper hold for tensors $\Cbb$ and ${\mathbb S}$ which may depend on $\Bx$ 
(i.e., which have spatial variations), and which may be anisotropic. 

Having established the equivalence of the boundary data, we can solve some inverse problems for the thin plate via the corresponding results for 2D elasticity, and of course, vice versa. For example, we can derive bounds on the volume fraction of the 2-phase thin plate via the result for 2D elasticity obtained by Milton and Nguyen \cite{ml}. On the other hand, the estimate of the size of an inclusion for the thin plate studied by Morrassi, Rosset, and Vessella \cite{mrv07, mrv09} is equivalent to the same problem for 2D elasticity, which was solved by Alessandrini, Morassi, and Rosset \cite{amr02} (corrected in \cite{amrs}). Moreover, we can also study the size estimate problem for 2D elasticity with certain anisotropic media through the similar result for the thin plate obtained by Morassi, Rosset, and Vessella in \cite{mrv11}.

The paper is organized as follows. In Section~2, we prove the explicit relations between the boundary data of the thin plate and that of 2D elasticity. In Section~3, we discuss the application of this equivalence to the size estimate problem for the thin plate equation.

\section{Equivalence of boundary data for the plate and for 2D elasticity}

It is well known that when the domain is simply connected, the plate
equation is equivalent to 2D elasticity equation, and vice versa. To
describe the equivalence,  let us define the rotational,
self-adjoint, 4th order tensor $\Rbb$ by
 \beq
 \Rbb \BA = \BR_{\perp}^T \BA \BR_{\perp}
 \eeq
for every $2\times 2$ matrix $\BA$, where
\[
\BR_{\perp}=\begin{bmatrix}0&1\\-1&0\end{bmatrix}.
\]
Then the equivalence between \eqref{plate} and \eqref{elasticity} is given by
 \beq\label{equiv}
 \Cbb = \Rbb {\mathbb S} \Rbb
 \eeq
and
 \beq\label{sigma}
 \Gs = \BR_{\perp}^T (\nabla^2 u) \BR_{\perp} \ (= \Rbb (\nabla^2 u)),
 \eeq
where in the context of  2D elasticity $u$ is known as the Airy stress function.
As a consequence of these two relations, we have
 \beq\label{eps}
 \Cbb \nabla^2 u = \BR_{\perp}^T \eps \BR_{\perp}.
 \eeq

Recall that the Neumann data for the thin plate equation are given by
 \[
 ({\mathbb C}\nabla^2u){\bf n}\cdot{\Bn}=-M_n, \quad \text{div}({\mathbb C}\nabla^2u)\cdot{\Bn} + (({\mathbb C}\nabla^2u){\bf n}\cdot{\Bt})_{,t}=(M_{t})_{,t}.
 \]
A comment on $(M_{t})_{,t}$ is helpful. Since $\text{div}\text{div}({\mathbb C}\nabla^2u)=0$, there exists a potential $\psi$ such that
\[
\text{div}({\mathbb C}\nabla^2u)=(\psi_{,2},-\psi_{,1}),
\]
and hence
\beq\label{psit}
\text{div}({\mathbb C}\nabla^2u)\cdot{\bf n}=\nabla\psi\cdot{\bf t}=\psi_{,t}\quad \text{on } \partial\Omega
\eeq
since ${\bf t}=-\BR_{\perp}{\bf n}$. Integrating \eqref{psit} along $\partial\Omega$ from some $x_0\in\partial\Omega$ and choosing an appropriate $\psi(\Bx_0)$, we obtain
\begin{equation}\label{f6}
\psi+({\mathbb C}\nabla^2u){\bf n}\cdot{\bf t}=M_t\quad\text{on}\ \partial\Omega.
\end{equation}

\subsection{Traction vs Dirichlet data}

Here, following \cite[p.158]{gu}, we connect the traction data for 2D elasticity to the Dirichlet data of the thin plate. It follows from \eqref{sigma} that
 \beq \label{2.7}
 \BR_\perp^T \Gs \Bn = (\nabla^2 u) \Bt = \begin{bmatrix} \nabla u_{,1} \cdot \Bt \\ \nabla u_{,2} \cdot \Bt \end{bmatrix}.
 \eeq
Thus by integrating $\BR_\perp^T \Gs \Bn$ along $\p\GO$, we recover $\nabla u=[u_{,1},u_{,2}]^T$ (up to a constant) on $\partial\Omega$. So $u_{,n}$ and $u_{,t}$ are recovered. By integrating $u_{,t}$ along $\p\GO$ we recover $u$ on $\p\GO$. On the other hand, from $u$ and $u_{,n}$ we determine $\nabla u$ on $\partial\Omega$ and therefore $\sigma{\bf n}$ through \eqref{2.7}.

\subsection{Displacement vs Moments}

Now we turn to the relation between the displacement of 2D elasticity and the moments of the thin plate. Since $\eps=\frac{1}{2} (\nabla \Bv + (\nabla\Bv)^T)$, we have
 \[
 \BR_{\perp}^T \eps \BR_{\perp} =
 \begin{bmatrix} v_{2,2} & -\frac{1}{2} v_{1,2} - \frac{1}{2} v_{2,1} \\
 -\frac{1}{2} v_{1,2} - \frac{1}{2} v_{2,1} & v_{1,1} \end{bmatrix}.
 \]
Thus we obtain
 \begin{align*}
 \div \, \BR_{\perp}^T \eps \BR_{\perp} & =
 \begin{bmatrix} v_{2,12} -\frac{1}{2} v_{1,22} - \frac{1}{2} v_{2,12} \\
 \nm
 v_{1,12} -\frac{1}{2} v_{1,21} - \frac{1}{2} v_{2,11} \end{bmatrix}
 = \begin{bmatrix} \frac{1}{2} v_{2,12} -\frac{1}{2} v_{1,22} \\
 \nm
 \frac{1}{2} v_{1,12} - \frac{1}{2} v_{2,11} \end{bmatrix} \\
 &= \BR_{\perp}^T \begin{bmatrix} (\frac{1}{2} v_{1,2} -\frac{1}{2} v_{2,1})_{,1} \\
 \nm
 (\frac{1}{2} v_{1,2} - \frac{1}{2} v_{2,1})_{,2} \end{bmatrix}
 = \BR_{\perp}^T \nabla (\frac{1}{2} v_{1,2} - \frac{1}{2} v_{2,1}),
 \end{align*}
which implies
 \begin{align*}
 \Bn \cdot (\div \, \BR_{\perp}^T \eps \BR_{\perp})
 &= (\BR_{\perp} \Bn) \cdot \nabla (\frac{1}{2} v_{1,2} - \frac{1}{2} v_{2,1}) \\
 &= - \Bt \cdot \nabla (\frac{1}{2} v_{1,2} - \frac{1}{2} v_{2,1}) \\
 &= (\frac{1}{2} v_{2,1} - \frac{1}{2} v_{1,2} )_{,t}.
 \end{align*}
It then follows from \eqref{psit} that
 $$
 \psi_{,t}= \Bn \cdot \text{div}({\mathbb C}\nabla^2u) = (\frac{1}{2} v_{2,1} - \frac{1}{2} v_{1,2} )_{,t},
 $$
and from \eqref{f6} that
 \beq\label{2.8}
 M_t = \frac{1}{2} (v_{2,1} - v_{1,2}) + (\BR_{\perp}^T \eps \BR_{\perp} \Bn) \cdot \Bt.
 \eeq
Observe that $\frac{1}{2} (v_{2,1} - v_{1,2})$ can be expressed as
 $$
 \frac{1}{2} (v_{2,1} - v_{1,2}) = \left(\BR_\perp^T \begin{bmatrix} 0 & -\frac{1}{2} (v_{2,1} - v_{1,2}) \\
 \frac{1}{2} (v_{2,1} - v_{1,2}) & 0 \end{bmatrix} \BR_\perp \Bn \right) \cdot \Bt.
 $$
We also have
 \begin{equation*}
 (\BR_{\perp}^T \eps \BR_{\perp} \Bn) \cdot \Bt = \left(\BR_\perp^T \begin{bmatrix} v_{1,1} & \frac{1}{2} v_{1,2} + \frac{1}{2} v_{2,1} \\
 \frac{1}{2} v_{1,2} + \frac{1}{2} v_{2,1} & v_{2,2} \end{bmatrix} \BR_\perp \Bn \right) \cdot \Bt.
 \end{equation*}
From \eqref{2.8} we thus find
 \begin{align}
 M_t &=  \left(\BR_\perp^T \begin{bmatrix} v_{1,1} & v_{1,2}\\
 v_{2,1} & v_{2,2} \end{bmatrix} \BR_\perp \Bn \right) \cdot \Bt \nonumber \\
 &= - \Bt \cdot \BR_\perp^T (\nabla \Bv) \Bt = - (\BR_\perp \Bt) \cdot (\nabla \Bv) \Bt = - \Bn \cdot (\nabla \Bv) \Bt. \label{Mt}
 \end{align}

We also deduce from the definition of $M_n$ that
 \beq\label{Mn}
 M_n = \Bn \cdot (\BR_{\perp}^T \eps \BR_{\perp}) \Bn = \Bt \cdot \eps \Bt = \Bt \cdot (\nabla \Bv) \Bt.
 \eeq
 Formulae \eqref{Mt} and \eqref{Mn} show that $M_n$ and $M_t$ can be recovered from $\Bv_{,t}$. On the other hand, we can recover $\Bv_{,t}$ from $M_t$ and $M_n$, i.e.,
$
\Bv_{,t}=-M_t{\bf n}+M_n{\bf t}.
$
 By integrating $\Bv_{,t}$ along $\p \GO$, we  then determine $\Bv$ on $\p\GO$.

\section{Applications to inverse problems}

\subsection{Estimating the volume of an inclusion for the thin plate with isotropic phases}

In this section, we want to discuss the size estimate problem for the thin plate equation based on the similar problem for 2D elasticity. Now we assume that the thin plate is made of an isotropic medium, i.e, the fourth order tensor ${\mathbb C}$ is given by
 \begin{equation}\label{C}
 C_{ijkl}=\frac{1}{2}B(1-\nu) (\delta_{ik}\delta_{jl}+\delta_{il}\delta_{jk})+B\nu\delta_{ij}\delta_{kl}.
 \end{equation}
The scalar $B$ is called the bending stiffness and is defined by
 \[
 B=\frac{h^3}{12}\left(\frac{E}{1-\nu^2}\right),
 \]
where $E$ is the Young's modulus, and $\nu$ is the
Poisson's coefficient. Both $E$ and $\nu$ can be written in terms of the
Lam\'e coefficients as follows:
$$E=\frac{\mu(2\mu+3\lambda)}{\mu+\lambda}\quad\text{and}\quad\nu=\frac{\lambda}{2(\mu+\lambda)}.$$
In terms of these the strong convexity condition \eqref{con} reads as
\[
\mu>\gamma\quad\text{and}\quad 2\mu+3\lambda>\gamma.
\]

Now assume that the plate is made by two different materials
in the sense that
$$\lambda=\lambda_1\chi_1+\lambda_2\chi_{2},\quad\mu=\mu_1\chi_1+\mu_2\chi_{2},$$ where
$$\chi_j=\begin{cases}
1\quad\text{in phase}\ j,\\
0\quad\text{otherwise}.
\end{cases}
$$
Suppose that $u$ is the solution of the thin plate equation \eqref{plate} having boundary data $\{u,u_{,n},M_n,M_t\}$ with
\[
 ({\mathbb C}\nabla^2u){\bf n}\cdot{\Bn}=-M_n, \quad \text{div}({\mathbb C}\nabla^2u)\cdot{\Bn} + (({\mathbb C}\nabla^2u){\bf n}\cdot{\Bt})_{,t}=(M_{t})_{,t}.
 \]
The inverse problem here is to estimate the volume fraction $f_1$ of phase 1 (or $f_2$) using $\{u,u_{,n},M_n,M_t\}$.

By \eqref{equiv}, we find that the thin plate equation with  elastic tensor \eqref{C}  is equivalent to 2D elasticity with a compliance tensor ${\mathbb S}=(S_{ijkl})$, where
\[
S_{ijkl}=\frac{1}{4\mu'}(\delta_{ik}\delta_{jl}+\delta_{il}\delta_{jk})+\frac 14(\frac{1}{\kappa}-\frac{1}{\mu'})\delta_{ij}\delta_{kl}
\]
where
\[
\kappa=\frac{2}{B(1+\nu)}\ (\text{bulk modulus}),\quad\mu'=\frac{2}{B(1-\nu)}\ (\text{shear modulus}).
\]
In other words, Hooke's law is given by
\[
\eps=\frac{1}{2\mu'}\sigma+\frac 14(\frac{1}{\kappa}-\frac{1}{\mu'})\text{Tr}(\sigma)\BI_2
\]
where $\BI_2$ is the identity matrix, and the equilibrium equation is
\[
\text{div}\ \sigma=0\quad\text{in}\quad\Omega.
\]
Using the result in the previous section, we can determine the displacement $\Bv$ and the traction $\sigma{\bf n}$ on $\partial\Omega$ from $\{u,u_{,n},M_n,M_t\}$. So to solve the size estimate problem for the thin plate, it suffices to solve the same problem for 2D elasticity.

Therefore, when all moduli involved are known, the result obtained by Milton and Nguyen \cite{ml} will give us bounds on the volume fraction and also the attainability conditions for these bounds for the two-phase thin plate with homogeneous isotropic phases. There are two approaches used in \cite{ml} -- the method of translation and the method of splitting. These methods depend on suitable null-Lagragians, which are (possibly nonlinear) functionals of fields that can be expressed in terms of the boundary measurements. For 2D elasticity, it is known that $\la \sigma\ra$, $\la\eps\ra$, $\la\text{det}\sigma\ra$, are null-Lagrangians, where $\la f\ra$ is the average of the field $f$. Thus, from relations \eqref{sigma} and \eqref{eps}, we immediately get that $\la\nabla^2u\ra$, $\la{\mathbb C}\nabla^2u\ra$, and $\la\text{det}\nabla^2u\ra$ are null-Lagrangians for the thin plate. It is interesting to point out that $\la\text{det}{\mathbb C}\nabla^2u\ra$ is not a null-Lagrangian for the thin plate. This corresponds to the fact that $\la\text{det}\eps\ra$ is not a null-Lagrangian for 2D elasticity.

We now consider the case where the two phases of the plate are themselves inhomogeneous, with spatially varying moduli. Assume that the medium in, say, phase 1, is given, but the medium in phase 2 is unknown. We also want to find upper and lower bounds on the volume of the unknown inclusion (phase 2) from one set of measurements of $\{u,u_{,n},M_n,M_t\}$. Under suitable conditions, this problem was solved in \cite{mrv07} (with the fatness assumption) and in \cite{mrv09} (for a general inclusion). The method used there is based on some quantitative estimates of the unique continuation property for the thin plate. Using the equivalence of inverse problems for the thin plate and for  2D elasticity, we can convert the size estimate problem for the
plate just described to the same problem for 2D elasticity, which was solved in \cite{amr02} (and also in \cite{amrs}).

\subsection{Estimating the volume of inclusion for 2D elasticity with anisotropic medium}

Now we assume that 2D elastic body is made of inhomogeneous and
\emph{anisotropic} medium with compliance tensor ${\mathbb
S}=(S_{ijkl})$. Due to the symmetry property for the compliance
tensor ${\mathbb S}$, we can denote
\[
\begin{cases}
S_{1111}=F,\ S_{1122}=S_{2211}=B,\\
S_{1112}=S_{1121}=S_{1211}=S_{2111}=-D,\\
S_{2212}=S_{2221}=S_{1222}=S_{2122}=-C,\\
S_{1212}=S_{1221}=S_{2112}=S_{2121}=E,\\
S_{2222}=A.
\end{cases}
\]
Suppose that the elastic body consists of two different media, i.e., ${\mathbb S}={\mathbb S}_0\chi_1+\widetilde{\mathbb S}\chi_2$. Likewise, the stress tensor $\sigma$ satisfies
\[
\text{div}\sigma=0\quad\text{and}\quad\eps={\mathbb S}\sigma\quad\text{in}\quad\Omega.
\]
Now assume that ${\mathbb S}_0$ is given and $\widetilde{\mathbb S}$ is unknown. We are now interested in estimating the size of $\chi_2$ by one boundary measurement $\{\Bv,\sigma{\bf n}\}$.

In view of the relation \eqref{equiv}, the corresponding elastic tensor ${\mathbb C}={\mathbb C}_0\chi_1+\widetilde{\mathbb C}\chi_2=(C_{ijkl})$ of the thin plate is given as
\[
\begin{cases}
C_{1111}=S_{2222}=A,\\
C_{1122}=S_{2211}=B,\\
C_{1112}=-S_{2212}=C,\\
C_{2212}=-S_{1112}=D,\\
C_{1212}=S_{1212}=E,\\
C_{2222}=S_{1111}=F.
\end{cases}
\]
So the problem above is equivalent to the size estimate problem for the thin plate with elastic tensor ${\mathbb C}$ using $\{u,u_{,n}.M_n,M_t\}$. This problem has been studied by Morassi, Rosset, and Vessella \cite{mrv11}. A key ingredient in the proof of \cite{mrv11} is the three-ball inequality for the plate equation with elastic tensor ${\mathbb C}_0$.  The three-ball inequality is proved under a so-called \emph{Dichotomy condition} for ${\mathbb C}_0$, which permits us to decompose the fourth order elliptic operator associated with the thin plate equation into a product of two second order elliptic operators. We now write this Dichotomy condition in terms of ${\mathbb S}_0=(S_{ijkl}^0)$.  We define
\[
\begin{cases}
S^0_{1111}=F_0,\ S^0_{1122}=S^0_{2211}=B_0,\\
S^0_{1112}=S^0_{1121}=S^0_{1211}=S^0_{2111}=-D_0,\\
S^0_{2212}=S^0_{2221}=S^0_{1222}=S^0_{2122}=-C_0,\\
S^0_{1212}=S^0_{1221}=S^0_{2112}=S^0_{2121}=E_0,\\
S^0_{2222}=A_0,
\end{cases}
\]
\[
a_0 =A_0,\quad a_1 =4C_0,\quad  a_2 =2B_0 +4E_0,\quad  a_3 =4D_0,\quad  a_4 =F_0,
\]
and the matrix
\[
{\mathcal M}(\Bx)=\begin{bmatrix}a_0&a_1&a_2&a_3&a_4&0&0\\0&a_0&a_1&a_2&a_3&a_4&0\\0&0&a_0&a_1&a_2&a_3&a_4\\4a_0&3a_1&2a_2&a_3&0&0&0\\0&4a_0&3a_1&2a_2&a_3&0&0\\0&0&4a_0&3a_1&2a_2&a_3&0\\0&0&0&4a_0&3a_1&2a_2&a_3\end{bmatrix}
\]
(Here we use similar notations as in \cite{mrv11}).  Then the Dichotomy condition is defined as
\begin{equation}\label{dic}
\frac{1}{a_0}\left|\text{det} {\mathcal M}(\Bx)\right|>0\quad\text{or}\quad \frac{1}{a_0}\left|\text{det} {\mathcal M}(\Bx)\right|=0\quad\forall\ \Bx\in\bar{\Omega}.
\end{equation}
Assume that \eqref{dic} holds and ${\mathbb S}_0\in C^{1,1}(\Omega)$ (equivalently, ${\mathbb C}^0\in C^{1,1}(\Omega)$). Then one can estimate the size of $\chi_2$ by one boundary measurement $\{\Bv,\sigma{\bf n}\}$ under the fatness assumption (see \cite[Theorem~3.2]{mrv11} for details).

\section*{Acknowledgements}
HK was partially supported by National Research Foundation of Korea
through NRF grants No. 2009-0085987 and 2010-0017532. GWM was
partially supported by the National Science Foundation of the USA
through grant DMS-0707978. JNW was partially supported by the
National Science Council of Taiwan through grants 100-2628-M-002-017
and 99-2115-M-002-006-MY3.

\end{document}